\documentclass[dvips]{imsart}
\usepackage{amsfonts}
\RequirePackage{xspace, ifthen}

\RequirePackage[OT1]{fontenc}
\RequirePackage{amsthm,amsmath,natbib}
\RequirePackage[colorlinks]{hyperref} \RequirePackage{hypernat}

\startlocaldefs \numberwithin{equation}{section}
\theoremstyle{plain}

\newtheorem{theorem}{Theorem}[section]

\newtheorem{lemma}[theorem]{Lemma}

\newtheorem{proposition}[theorem]{Proposition}
\newtheorem{corollary}[theorem]{Corollary}

\newcommand{\convl}[1]{\stackrel{#1}{\longrightarrow}}
\def\cip{\convl{\rm p   }{}}

\def\ind{\mbox{\bf\large $1\hspace{-0.3em}{\rm I}$}\xspace}

\def\scd{\ensuremath{{\mathcal{D}}}\xspace}
\def\scf{\ensuremath{{\mathcal{F}}}\xspace}

\def\scp{\ensuremath{{\mathcal{P}}}\xspace}
\def\scm{\ensuremath{{\mathcal{M}}}\xspace}
\def\scy{\ensuremath{{\mathcal{Y}}}\xspace}
\def\Del{\Delta}
\def\summ#1#2#3{\sum_{#1=#2}^{#3}}
\def\argmin{\mathop{\mathrm{arg\,min}}}
\def\E{{\rm E}\mathop{\!}\nolimits}
\def\dist{\sim}
\newcommand{\eqsplit}[2][*]%
  {\ifthenelse{\equal{#1}{*} \or\equal{#1}{<???>}
                     \or\equal{#1}{nn} \or\equal{#1}{ } \or\equal{#1}{}}
    { \begin{align*}%
         #2 %
         \end{align*}%
        }
    {\begin{equation}\label{#1}\begin{split}\allowdisplaybreaks%
         #2%
         \end{split}\end{equation}
        }
  }

\def\iid{i.i.d.\xspace }
\def\bb#1{\mathbb{#1}}
\def\bbg{\ensuremath{{\bb G}}\xspace}
\def\bbk{\ensuremath{{\bb K}}\xspace}
\def\eps{\varepsilon}
\def\ti{\tilde}
\def\fun#1{\mathop{\rm #1}\nolimits}
\def\var{\fun{Var}}
\def\en{\ensuremath{\infty}\xspace}

\endlocaldefs

\begin{document}

\begin{frontmatter}
\title{Asymptotic efficiency of simple decisions for the compound decision problem\protect\thanksref{T1}}
\runtitle{simple estimators for compound decisions}
\thankstext{T1}{Partially supported by NSF grant DMS-0605236, and an ISF Grant. }


\begin{aug}
\author{\fnms{Eitan.} \snm{ Greenshtein}
\ead[label=e1]{eitan.greenshtein@gmail.com}} \and
\author{\fnms{Ya'acov} \snm{Ritov}
\ead[label=e2]{yaacov.ritov@gmail.com}}

\runauthor{Greenshtein and Ritov}

\affiliation{}

\address{Department of Statistical Sciences\\Duke University\\ Durham,  NC 27708-0251, USA\\
\printead{e1} \phantom{E-mail:\ }}

\address{Jerusalem, Israel\\
\printead{e2} \phantom{E-mail:\ }}

\end{aug}

\begin{abstract}
We consider the compound decision problem of estimating
a vector of $n$ parameters, known up to a permutation, corresponding to $n$ independent observations,
and discuss the difference between two symmetric classes of estimators.
The first and larger class is restricted to the set of all permutation invariant estimators. The second class is restricted further to simple symmetric procedures. That is, estimators such that each parameter is estimated by a function of the corresponding observation alone. We show that under mild conditions, the minimal total squared error risks over these two classes are asymptotically equivalent up to essentially $O(1)$ difference.

\end{abstract}



\end{frontmatter}


\def\bbmu{\setbox\ewa\hbox{\mu}\setlength{\deno}{\wd\ewa} \mu\hspace{-0.95\deno}\mu}
\def\bbmus{\setbox\ewa\hbox{\mu}\setlength{\deno}{\wd\ewa} \mu\hspace{-0.75\deno}\mu}
\def\bmu{\ensuremath{\boldsymbol{\mu}}\xspace}

\section{Introduction{}}

Let $\scf=\{F_\mu: \mu \in \scm \}$ be a parametrized family of distributions.
Let $Y_1,Y_2\dots$ be a sequence of independent random variables,
where $Y_i\in\scy$ and $Y_i\dist F_{\mu_i}$, $i=1,2,\dots $.
For each $n$, we suppose that the sequence $\mu_{1:n}$ is known up to a permutation, where for any sequence $x=(x_1,x_2,\dots)$ we denote the sub-sequence $x_s,\dots,x_t$ by $x_{s:t}$. We denote by $\bmu=\bmu_n$ the set $\{\mu_1,\dots,\mu_n\}$, i.e., $\bmu$ is $\mu_{1:n}$ without any order information. We
consider in this note the problem of estimating it by
$\hat{\mu}_{1:n}$ under the loss $ \summ i1n (\hat{\mu}_i - \mu_i)^2$,
where   $\hat{\mu}_{1:n}= \Delta(Y_{1:n})$.
We assume that the family $\scf$ is dominated by a measure $\nu$, and denote the corresponding densities simply by $f_{i}=f_{\mu_i}$, $i=1,...,n$.
The important example is, of course,   $F_{\mu_i}=N(\mu_i,1)$.

Let $\scd^S=\scd^S_{n}$ be the set of all \emph{simple symmetric decision functions} $\Delta$, that is, all $\Delta$  such that  $\Delta(Y_{1:n})= (\delta(Y_1),\dots,\delta(Y_n))$, for some function $\delta:\scy\to\scm$.
In particular, the best simple symmetric function is denoted by $\Del^S_{\bmu}=(\delta^S_{\bmu}(Y_1),\dots,\delta^S_{\bmu}(Y_n))$:
$$\Delta^S_{\bmu}= \argmin_{\Delta \in \scd^S_n} \E || \Delta-\mu_{1:n}||^2,$$
and denote $$r^S_n= \E || \Delta^S_{\bmu}(Y_{1:m})-\mu_{1:n}||^2, $$
where, as usual, $\|a_{1:n}\|^2=\summ i1n a_i^2$.

The class of simple rules may be considered too restrictive. Since the
$\mu$s are known up to a permutation, the problem seems to be of
matching the $Y$s to the $\mu$s. Thus, if $Y_i\dist N(\mu_i,1)$, and
$n=2$, a reasonable decision would make
$\hat\mu_1$ closer to $\mu_1\wedge\mu_2$ as $Y_2$ gets
larger. The simple rule clearly remains inefficient if the $\mu$s are well separated, and generally speaking, a bigger class of decision rules may be needed to obtain efficiency. However, given the natural invariance of the problem, it makes sense to be restricted to the class  $\scd^{PI}=\scd^{PI}_n$  of all permutation invariant decision functions, i.e, functions $\Delta$
that satisfy  for any permutation $\pi$ and any $(Y_1,...,Y_n)$:
\eqsplit{
\Delta(Y_{1},\dots,Y_n)=(\hat{\mu}_1,...,\hat{\mu}_n) \quad\iff\quad \Delta(Y_{\pi(1)},..,Y_{\pi(n)} )=(\hat{\mu}_{\pi(1)},...,\hat{\mu}_{\pi(n)} ).
}

Let $$\Delta^{PI}_{\bmu}= \argmin_{ \Delta \in \scd^{PI}} \E ||\Delta(Y^n)-\mu_{1:n}||^2$$
be the optimal permutation invariant rule under $\bmu$, and denote its risk by
$$r^{PI}_n= E || \Delta^{PI}_{\bmu} (Y_{1:n}) - \mu_{1:n}||^2.$$

Obviously ${\cal D}^S \subset {\cal D}^{PI}$, and whence $r_n^S \geq r_n^{PI}$.
Still, `folklore', theorems in the spirit of De Finetti, and results
like  Hannan and Robbins (1955),
imply that asymptotically (as $n \rightarrow \infty$)
$\Delta^{PI}_{\mu^n}$ and $\Delta^S_{\mu^n}$ will have `similar' risks.
Our main result establishes conditions that imply $$r^S_n- r^{PI}_n =O(1).$$
To repeat, \bmu is assumed known in this note. In the general decision theory framework the unknown parameter is the order of its member to correspond with $Y_{1:n}$, and the parameter space, therefore, corresponds to the set of all the permutation of $1,\dots,n$.

An asymptotic equivalence as above implies, that when we confine ourselves to the class of permutation invariant procedures, we may further restrict ourselves to the class of  simple symmetric procedures, as is usually done in  the standard analysis of compound decision problems. The later class is smaller and simpler.

The motivation for this paper stems from the way the notion of oracle is used in some sparse estimation problems.  Consider two oracles both {\it know} the value of $\bmu$. Oracle I is restricted to use only a procedure from the class $ {\cal D}^{PI}$, while Oracle II is restricted to use only procedures from ${\cal D}^{S}$. Obviously Oracle I has an advantage, our results quantify this advantage and show that it is asymptotically negligible. Furthermore, starting with Robbins (1951) various
oracle-inequalities were obtained showing that one can achieve
nearly the risk of Oracle II, by a `legitimate' statistical procedure.
See, e.g., the survey Zhang (2003), for  oracle-inequalities regarding
the difference in risks. See also Brown and Greenshtein (2007), and
Jiang and Zhang (2007) for oracle inequalities regarding the ratio
of the risks. However, Oracle II is weak, and hence, these claims may seem to be too weak. Our equivalence results, extend many of those oracle inequalities to be valid also with respect to Oracle I.  We needed a stronger result than the usual objective that the mean risks are equal up to $o(1)$ difference. Many of the above mentioned recent applications of the compound decision notion are about sparse situations when most of the $\mu$s are in fact 0,  the mean risk is $o(1)$, and the only interest is in total risk.

Let $\bmu_1,\dots,\bmu_n$ be some arbitrary ordering of $\bmu$. Consider now the Bayesian model under which $(\pi,Y_{1:n})$, $\pi$ a random permutation, have a distribution given by
 \eqsplit[Bayesian] {
    &\text{$\pi$ is uniformly distributed over $\scp(1:n)$;}
    \\
    &\text{Given $\pi$, $Y_{1:n}$ are independent, $Y_i\dist F_{\bmu_{\pi(i)}}$, $i=1,\dots,n$,}
    }
where for every $s<t$, $\scp(s:t)$ is the set of all permutations of $s,\dots,t$.  The above description induces a joint distribution of $(M_1,...,M_n,Y_1,...,Y_n)$, where  $ M_i \equiv \bmu_{\pi(i)},$
for a random permutation $\pi$.

The first part of the following proposition is a simple special case
of general theorems representing the best invariant procedure under
certain groups, as Bayes with respect to the appropriate Haar
measure; for background see, e.g.,   Berger (1985), Chapter 6.
The second part of the proposition was derived in various papers starting with Robbins (1951).

In the following proposition and proof, $\E_{\mu_{1:n}}$ is the expectation under the model in which the observations are independent and $Y_i\dist F_{\mu_i}$ $\E_{\bmu}$ is the expectation under the above joint distribution  of $Y_{1:n}$ and $M_{1:n}$. Note that under the latter model, marginally $M_i\sim\bbg_n$, the empirical measure defined by \bmu, and conditional on $M_i=m$, $Y_1\dist F_{m}$, $i=1,\dots,n$.
\begin{proposition}
\label{prop:basic}
The best simple and permutation invariant rules are given by
\begin{enumerate}
\item[(i)]
$\Delta^{PI}_{\bmu}(Y_{1:n})= E_{\bmu}\bigl( M_{1:n}\bigr| Y_{1:n} \bigr)$.
\item[(ii)]
$\delta^S_{\bmu}(Y_i)= E_{\bmu}(M_i|Y_i)$, $i=1,\dots,n$.
\item[(iii)]
Suppose  $\E_{\bmu}\|\Del^{S}_{\bmu} - \Del^{PI}_{\bmu}\|^2=O(\eps_n^2)$, then $r_n^{S}=r_n^{PI}+O(\eps_n^2)$.
\end{enumerate}
\end{proposition}
\begin{proof}
We need only to  give the standard proof of the third part.
First, note that by invariance $\Delta_{\bmu}^{PI}$ is an equalizer
(over all the permutations of $\bmu$), and hence  $\E_{\mu_{1:n}}(\Delta_{\bmu}^{PI}-
\mu_{1:n})^2= \E_{\bmu}(\Delta_{\bmu}^{PI}-M_{1:n})^2$. Also
$\E_{\mu_{1:n}}(\Delta_{\bmu}^S- \mu_{1:n})^2= \E_{\bmu}(\Delta_{\bmu}^S- M_{1:n})^2$.
Then, given the above joint distribution,
 \eqsplit{
    r_n^S &=   \E_{\bmu} \|\Del^S_{\bmu} - M_{1:n}\|^2
    \\
    &=  \E_{\bmu} \E_{\bmu}\{\|\Del^S_{\bmu} - M_{1:n}\|^2|Y_{1:n}\}
    \\
    &=   \E_{\bmu} \E_{\bmu}\{\|\Del^S_{\bmu} - \Del^{PI}_{\bmu}\|^2 + \|\Del^{PI}_{\bmu} - M_{1:n}\|^2|Y_{1:n}\}
    \\
    &= r_n^{PI} + O(\eps_n^2).
  }
\end{proof}

We now briefly review some related literature and problems. On simple
symmetric functions, compound decision and its relation to empirical
Bayes,
see Samuel (1965), Copas (1969), Robbins (1983), Zhang (2003), among many other
papers.

Hannan and Robbins (1955) formulated essentially the same
equivalence problem in testing problems, see their Section 6. They
show for a special case an equivalence up to
 $o(n)$ difference in the `total risk' (i.e., non-averaged risk).
Our results for estimation under squared loss are stated in terms of
the total risk and we obtain $O(1)$ difference.

Our results have a strong connection to De Finetti's Theorem. The exchangeability
induced on $M_1,...,M_n$, by the Haar measure, implies `asymptotic
independence' as in Definetti's theorem, and consequently asymptotic
independence of $Y_1,...,Y_n$. Thus we expect $E(M_1|Y_1)$ to be
asymptotically similar to $E(M_1|Y_1,...,Y_n)$. Quantifying this
similarity as $n$ grows,
has to do with the rate of convergence in DeFinetti's theorem. Such
rates were established by Diaconis and Freedman (1980), but are not
directly applicable to obtain our results.

After quoting a simple result in the following section, we consider in  Section \ref{sec:bin} the special important, but simple, case of two-valued parameter. In Section \ref{sec:cont} we obtain a strong result under strong
conditions.  Finally, the main result is given in Section \ref{sec:main},   it covers the two preceding cases, but with some price to pay for the generality.

\section{Basic lemma and notation}
\label{sec:basic}
The following lemma is standard in comparison of experiments theory;
for background on comparison of experiments in testing see Lehmann
(1986), p-86.
The proof follows a simple application of Jensen's inequality.
\begin{lemma}
\label{lem:ExpComp}
Consider two pairs of distribution, $\{G_0,G_1\}$ and  $\{\ti G_0, \ti G_1\}$. Suppose that there exists a Markov kernel \bbk such that $G_i(\cdot)=\int \bbk(y,\cdot)\,d\ti G_i(y)$, $i=1,2$. Then
 \eqsplit{
    \E_{G_0} \psi\bigl( \frac {dG_1}{dG_0}\bigr) & \leq \E_{\ti G_0}
    \psi\bigl( \frac { d \tilde{G}_1}{ d \tilde{G}_0}\bigr)
  } for any convex function $\psi$
\end{lemma}

For simplicity  denote $f_i(\cdot)=f_{\mu_i}(\cdot)$, and for any
random variable $X$, we may write $X\dist g$ if $g$ is its density
with respect to a certain dominating measure. Finally, for simplicity we use the notation $y_{-i}$ to denote the sequence $y_1,\dots,y_n$ without its $i$ member, and similarly $\bmu_{-i} =\{\mu_1,\dots,\mu_n\}/\mu_i$. Finally $f_{-i}(Y_{-j})$ is the marginal density of $Y_{-j}$ under the model \eqref{Bayesian} conditional on $M_j=\mu_i$.


\section{Two valued parameter}
\label{sec:bin}

We suppose in this section that $\mu$ can get one of two values which we denote by $\{0,1\}$.  To simplify notation we denote the two densities by $ f_0$ and $ f_1$.
\begin{theorem}
Suppose that either of the following two conditions hold: 
\begin{enumerate}
\item[(i)]
$f_{1-\mu}(Y_1)/f_{\mu}(Y_1)$ has a finite variance under $\mu\in\{0,1\}$.
\item[(ii)]
$\summ i1n \mu_i/n\to\gamma\in(0,1)$, and
$f_{1-\mu}(Y_1)/f_{\mu}(Y_1)$ has a finite variance under $\mu=0$  or $\mu=1$.
\end{enumerate}
Then $\E\|\hat\mu^S-\hat\mu^{PI}\|^2=O(1)$.
\end{theorem}

\begin{proof}
Suppose condition (i) holds.
Let $K=\summ i1n \mu_i$, and suppose, wlog, that $K\leq n/2$. Consider the Bayes model of \eqref{Bayesian}. By Bayes Theorem
 \eqsplit{
    P(M_1=1|Y_{1} )
    &= \frac { K  f_1(Y_1) }{   K  f_1(Y_1)   + (n-K)  f_0(Y_1) } ,
    }
where, with some abuse of notation $f_k(Y_{2:k})$ is the joint density of $Y_{2:n}$ conditional on $\summ j2n \mu_j=k$.
On the other hand
    \eqsplit{
    &P(M_1=1|Y_{1:n} )\\
    &= \frac {K  f_1(Y_1) f_{K-1}(Y_{2:n})  }{ K  f_1(Y_1) f_{K-1}(Y_{2:n})  + (n-K)  f_0(Y_1) f_K(Y_{2:n}) }
    \\
    &= \frac { K  f_1(Y_1) }{   K  f_1(Y_1)   + (n-K)  f_0(Y_1) }
    \Bigl(1+  \frac { (n-K) f_0(Y_1) }{   K  f_1(Y_1)   + (n-K)  f_0(Y_1) } \bigl(\frac{f_{K}}{f_{K-1}}(Y_{2:n})-1\bigr)    \Bigr)^{-1}
    \\
    &=  P(\mu_{\pi(1)}=1|Y_{1} )         \Bigl(1+O\bigl(\frac{f_{K}}{f_{K-1}}(Y_{2:n})-1\bigr)\Bigr)
  }

We use Lemma \ref{lem:ExpComp} to compare the testing between $f_K(Y_{2:k})$ vs. $f_{K-1}(Y_{2:k})$ to an easier problem, from which the original problem can be obtained by adding a random permutation. Suppose for simplicity and wlog that in fact $Y_{2:K}$ are \iid under $f_1$, while $Y_{K+1:n}$ are \iid under $f_0$. Then we compare
 \eqsplit{
    g_{K-1}(Y_{2:n}) &=  \prod_{j=2}^K f_1(Y_j) \prod _{j=K+1}^n f_0 (Y_j),
  }
the true distribution, to the mixture
 \eqsplit{
    g_K(Y_{2:n}) &= g_{K-1}(Y_{2:n}) \frac1{n-K} \summ j{K+1}n \frac{f_1}{f_0}(Y_j).
  }
However, the likelihood ratio between $g_K$ and $g_{K-1}$ is a sum of $n-K$ terms, each with mean 0 (under $g_{K-1}$) and finite variance. It is, therefore, $1+O_p(n^{-1/2})$ in the mean square.

Consider now the second condition. By assumption, $K$ is of the same order as $n$, and we can assume, wlog, that the $f_1/f_0$ has a finite variance under $f_0$. With this understanding, the above proof holds for the second condition.
\end{proof}

The condition of the theorem is clearly satisfied in the normal shift model: $F_i=N(\mu_i,1)$, $i=1,2$. It is satisfied for the normal scale model,  $F_i=N(0,\sigma^2_i)$, $i=1,2$,  if $K$  is of the same order as $n$, or if $\sigma_0^2/2 < \sigma_1^2 < 2\sigma_0^2$.


\section{Dense $\mu$s }
\label{sec:cont}

We consider now another simple case in which $\bmu$ can be ordered $\mu_{(1)},\dots,\mu_{(n)}$ such that the difference $\mu_{(i+1)}-\mu_{(i)}$ is uniformly small. This will happen if, for example, $\bmu$ is in fact a random sample from a distribution with density with respect to Lebesgue measure, which is bounded away from 0 on  its support, or more generally, if it sampled from a distribution with short tails. Denote by $Y_{(1)},\dots,Y_{(n)}$ and  $f_{(1)},\dots,f_{(n)}$ the corresponding ordering of the $Y$s and $f$s.

We assume in this section
\begin{enumerate}
\item[(B1)]For some slowly converging to infinity constants $A_n$ and $V_n$:
 \eqsplit{
    \max_{i,j} |\mu_i-\mu_j| &= A_n
    \\
    \var \Bigl(\frac{f_{(j+1)}}{f_{(j)}} (Y_{(j)})\Bigr)&\leq \frac{V_n}{n^2}.
  }
\end{enumerate}
Note that condition \textbf{(B1)} holds for both the normal shift model and the normal scale model, if \bmu behaves like a sample from a distribution with a density as above.
\begin{theorem}
If Assumption \textbf{(B1)} holds then
 \eqsplit{
    \summ i1n |\mu_i^{PI} - \mu_i^{S} |^2 = O_p(A_n^2V_n^2/n).
  }
\end{theorem}

\begin{proof}
By definition
 \eqsplit{
    \mu_1^{S} &= \frac{\summ i1n \mu_i f_i(Y_i)} {\summ i1n  f_i(Y_i)}
    \\
    \mu_1^{PI} &= \frac{\summ i1n \mu_i f_i(Y_i)f_{-i}(Y_{2:n})} {\summ i1n  f_i(Y_i)f_{-i}(Y_{2:n})} }
where $f_{-i}$ is the density of $Y_{2:n}$ under $\bmu_{-i}$:
 \eqsplit{
    f_{-i}((y_{2:m}))  &= \frac1{(n-1)!} \sum_{\pi\in\scp(2:n)} \prod_{j=2}^nf_{\pi(j)}(y_j) \frac{f_1}{f_{i}}(y_i).
    }

The result will follow if we argue that $\max_i|f_{-i}(Y_{2:n})/f_{-1}(Y_{2:n})-1|\cip 1$. In fact we will establish a slightly stronger claim that
 \eqsplit{
    \|f_{-i} - f_{-1}\|_{TV} \to 0
  }
where $\|\cdot\|_{TV}$  denotes the total variation norm.

We will bound this distance by the distance between two other densities. Let $g_{-1}(y_{2:n}) = \prod_{j=2}^n f_j(y_j) $, the true distribution of $Y_{2:n}$. We define now a similar analog of $f_{-i}$. Let $r_j$ and  $y_{(r_j)}$  be defined by $f_{j}=f_{(r_j)}$ and $y_{(r_j)}=y_j$, $j=1,\dots,n$. Suppose, for simplicity, that $r_i<r_1$. Let
 \eqsplit{
    g_{-i}(y_{2:n})&= g_{-1}(y_{2:n}) \prod_{j=r_i}^{r_1-1} \frac{f_{(j+1)}}{f_{(j)}} (y_{(j)}).
  }
Note that $g_{-i}$    depends only on $\bmu_{-i}$. Moreover, if $\ti Y_{2:n}\dist g_{-j}$, then one can obtain $Y_{2:n}\dist f_{-j}$ by the Markov kernel that takes $\ti Y_{2:n}$ to a random permutation of itself. It follows from Lemma \ref{lem:ExpComp}
 \eqsplit{
    \|f_{-i} - f_{-1}\|_{TV} &\leq \|g_{-i} - g_{-1}\|_{TV}
    \\
    &= \E_{\mu_{2:n}} \Bigl|  \frac{g_{-i}}{  g_{-1}} (Y_{2:n}) -1\Bigr|
    \\
    &= \E_{\mu_{2:n}} \Bigl|  \prod_{j=k}^{r_1-1}\frac{f_{(j+1)}}{f_{(j)}} (Y_{(j)})  -1\Bigr|
  }
But, by assumption
 \eqsplit{
    R_k &=    \prod_{j=k}^{r_1-1} \frac{f_{(j+1)}}{f_{(j)}} (Y_{(j)})
  }
is a reversed $L_2$ martingale, and it follows from Assumption \textbf{(B1)} that
 \eqsplit{
    \max_{k<r_1} |R_k -1|  &= O_p(A_nV_n/n) .
  }
Similar argument applies to $i$, $r_i>r_1$, yielding
 \eqsplit{
    \max_{i}\|f_{-i} - f_{-1}\|_{TV} = O_p(A_nV_n/n)
  }
 But then we argue
 \eqsplit{
    |\mu_1^{PI} - \mu_1^{S} | \leq \max_{i,j} |\mu_i-\mu_j| \bigl( \max_{i,j} \frac{f_{-i}} {f_{-j}} (Y_{2:n}) - 1\bigr) = O_p(A_nV_n/n).
  }
  The theorem follows.
\end{proof}

\section{Main result}
\label{sec:main}

We assume that for some $C<\en$:
\begin{enumerate}
\item[(G1)]
We have $\max_{i\in\{1,\dots,n\}}|\mu_i|<C $  and $\max_{i,j\in{1,\dots,n}}\E_{\mu_i}( f_{\mu_j}(Y_1) /f_{\mu_i}(Y_1) )^2<C$. Finally there is $\gamma>0$ such that $\min_{i,j\in{1,\dots,n}}P_{\mu_i}( f_{\mu_j}(Y_1) /f_{\mu_i}(Y_1) > \gamma )\geq 1/2$.
\item[(G2)] Let
 \eqsplit{
    p_j(Y_i) &=  \frac{ f_j(Y_i)}  {\summ k1n f_k(Y_i)},\quad i,j=1,\dots,n.
  }
 Then
 \eqsplit{
    \E \summ i1n \summ j1n \bigl(p_j(Y_i)\bigr)^2 &< C
    \\
    \summ i1n\E \frac 1{n \min_j p_j(Y_i)}    &<Cn
    \\
    \E\summ i1n \frac {\summ j1n \bigl(p_j(Y_i)\bigr)^2}{n \min_j p_j(Y_i)}    &<C.
  }

\end{enumerate}

Both assumptions describe a situation where the $\mu$s do not ``separate". They cannot be too far one from another, geometrically or statistically (Assumption \textbf{(G1)}), and they are dense in the sense that each $Y$ can be explained by many of the $\mu$s (Assumption \textbf{(G2)}). The conditions hold for the normal shift model if $\bmu_n$ are uniformly bounded: Suppose the common variance is 1 and $|\mu_j|<A_n$. Then
 \eqsplit{
    \E  \summ j1n \Bigl(\frac{ f_j(Y_1)}  {\summ k1n f_k(Y_1)}\Bigr)^2
    &= \E \frac{\summ j1n f_j^2(Y_1)} {(\summ k1n f_k(Y_1))^2}
    \\
    &\leq \E \frac{n e^{-Y_1^2+2A_nY_1-A_n^2}} {(n e^{-(Y_1^2-2A_nY_1+A_n^2)/2})^2}
    \\
    &= \frac1n \E e^{4A_nY_1}
    \\
    & = \frac1n e^{8A_n^2+4A_n\mu_1} \leq \frac1n e^{12A_n^2}.
  }
and the first part of \textbf{(G2)} hold. The other parts follow a similar calculations.
\begin{theorem}
Assume that \textbf{(G1)} and  \textbf{(G2)} hold. Then
 \eqsplit{
     \E\| \Del^{S}_{\bmu} - \Del^{PI}_{\bmu} \|^2 &= O(1) \tag{i}
     \\
     r^S_n-r^{PI}_n&= O(1). \tag{ii}
  }
\end{theorem}
\begin{corollary}
Suppose $\scf=\{N(\mu,1):\;|\mu|<c\}$ for some $c<\en$, then the conclusions of the theorem follow.
\end{corollary}
\begin{proof}
It was mentioned already in the introduction that when we are restricted to permutation invariant procedure we can consider the Bayesian model under which $(\pi,Y_{1:n})$, $\pi$ a random permutation, have a distribution given by \eqref{Bayesian}. Fix now $i\in\{1,\dots,n\}$. Under this model  we want to compare
 \eqsplit{
    \mu^S_i &= E(\mu_{\pi(i)}|Y_i),\quad i=1,\dots,n
  }
to
 \eqsplit{
    \mu^{PI}_i  &= E(\mu_{\pi(i)}|Y_{1:n}),\quad i=1,\dots,n .
  }
More explicitly:
 \eqsplit[mus] {
    \mu^S_i &= \frac{\summ j1n \mu_j f_j(Y_i)}{\summ j1n  f_j(Y_i)}
    \\
    &= \summ j1n \mu_j p_j(Y_i), \quad i=1,\dots,n
    \\
    \mu^{PI}_i  &= \frac{\summ j1n \mu_j f_j(Y_i) f_{-j} (Y_{-i}) } {\summ j1n  f_j(Y_i) f_{-j} (Y_{-i}) }
    \\
    &=  \summ j1n \mu_j p_j(Y_i) W_j(Y_{-i},Y_i ), \quad i=1,\dots,n,
  }
where  for all $i,j=1,\dots,n$, $f_j(Y_i)$ was defined in Section \ref{sec:basic}, and
 \eqsplit{
    p_j(Y_i) &= \frac{f_j(Y_i)} { \summ k1n f_k(Y_i)} ,
    \\
    W_j(Y_{-i},Y_i) &= \frac {f_{-j} (Y_{-i}) } {\summ k1n p_k(Y_i) f_{-k} (Y_{-i}) }
  }

Note  that $\summ k1n p_k(Y_i)=1$, and $W_j(Y_{-i},Y_i)$ is the likelihood ratio between two (conditional on $Y_i$) densities of $Y_{-i}$, say $g_{j0}$  and $g_{1}$.  Consider two other densities (again, conditional on $Y_i$):
 \eqsplit{
    \ti g_{j0}(Y_{-i}|Y_i) &=  {f_i}  (Y_j) \prod_{m\ne i,j}  f_m (Y_m),
    \\
    \ti g_{j1} (Y_{-i}|Y_i) &=  \ti g_{j0}(Y_{-i}|Y_i) \Bigl( \sum_{k\ne i,j} p_k(Y_i)\frac{f_i} {f_k} (Y_k)  + p_i(Y_i)\frac{f_j} {f_i} (Y_j) +p_j(Y_i) \Bigr)
  }
Note that $ g_{j0} = \ti g_{j0} \circ \bbk $ and  $ g_{1} = \ti g_{j1} \circ \bbk $, where $\bbk$ is the Markov kernel that takes $Y_{-i}$ to a random permutation of itself. It follows from Lemma \ref{lem:ExpComp} that
\eqsplit[ewi] {
    \E \bigl(|W_j(Y_{-i},Y_i)-1|^2\bigl|Y_i\bigr) &\leq \E_{\ti g_{j1}} \Bigl| \frac {\ti g_{j0}} {\ti g_{j1}}-1 \Bigr|^2
    \\
    &=  \E_{\ti g_{j0}} \Bigl| \frac {\ti g_{j0}} {\ti g_{j1}} -2 + \frac { \ti g_{j1}} {\ti g_{j0}} \Bigr|.
  }
This expectation does not depend on $j$: $i$ is related only to the $Y$s while $j$ is related only to $\bmu$.  Hence, to simplify notation, we take wlog $j=i$. Denote
 \eqsplit{
    L &= \frac { \ti g_{i1}} {\ti g_{i0}}
    = p_i(Y_i) +  \sum_{k\ne i} p_i(Y_i) \frac{f_i}{f_k}(Y_k)
    \\
    V&= \frac n4\gamma \min_jp_j(Y_i) ,
  }
where $\gamma$ is as in \textbf{(G1)}. Then by \eqref{ewi}
\eqsplit[ewiII]{
    \E \bigl(|W_j(Y_{-i},Y_i)-1|^2\bigl|Y_i\bigr) &\leq
    \E_{\ti g_{j0}} \Bigl(\frac1L-2+L \Bigr)
    \\
    &=  \E_{\ti g_{j0}} \frac{(L-1)^2}L
    \\
    & \leq  \frac1V \E_{\ti g_{j0}} (L-1)^2\ind(L>V) + \E_{\ti g_{j0}} \frac{\ind(L\leq V)}L
    \\
    &\leq     \E_{\ti g_{j0}} \frac{\ind(L\leq V)}L + \frac1{V}\summ k1n p_k^2(Y_i),
  }
by \textbf{G1}. Bound
 \eqsplit{
 L \geq  \gamma \min_k p_k(Y_i) \summ k1n  \ind(\frac{f_i}{f_k}(Y_k)>\gamma)\geq \gamma \min_k p_k(Y_i) (1+U),
      }
where  $U\dist B(n-1,1/2)$ (the 1 is for  the $i$th summand). Hence
 \eqsplit[ool]{
     \E_{\ti g_{j0}} \frac{\ind(L\leq V)}L
     &\leq \frac 1{\gamma \min_k p_k(Y_i)} \summ k0{\lceil n/4\rceil}
     \frac1{k+1} \binom{n-1}k  2^{-n+1}
     \\
     &=  \frac 1{\gamma n\min_k p_k(Y_i)}\summ k0{ \lceil n/4 \rceil} \binom{n}{k+1}2^{-n+1}
     \\
     &= O(e^{-n})  \frac 1{\gamma n\min_k p_k(Y_i)}
  }
by large deviation.

From \textbf{(G1)},  \textbf{(G2)}, \eqref{mus}, \eqref{ewiII}, and \eqref{ool}:
 \eqsplit{
    \E\E\bigl((\mu_i^S-\mu_i^{PI})^2\bigr|Y_i\bigr)
         &= \E\E\biggl(\Bigl(\summ j1n \mu_j p_j(Y_i)\bigl( W_j(Y_{-i},Y_i )-1\bigr)\Bigr)^2\bigr|Y_i\biggr)
         \\
         &\leq \max_j |\mu_j|^2 \E\summ j1n  p_j(Y_i)\E\bigl( W_j(Y_{-i},Y_i )-1\bigr)^2\Bigr)\bigr|Y_i\biggr)
         \\
     &\leq \kappa C^3/n,
  }
for some $\kappa$ large enough. Claim (i) of the theorem follows. Claim (ii) follows (i) by Proposition \ref{prop:basic}.

\end{proof}

\end{document}